\documentclass[twosided,reqno]{amsart}

\usepackage{amsmath}
\usepackage{amssymb}
\usepackage{amsthm}

\theoremstyle{theorem}
\newtheorem{theorem}{\scshape Theorem }[section]

\theoremstyle{definition}

\numberwithin{equation}{section}

\begin{document}

\title[On the degenerate Frobenius-Euler polynomials]{On the degenerate Frobenius-Euler polynomials}

\author{Taekyun Kim}
\address{Department of Mathematics, Kwangwoon University, Seoul 139-701, Republic of Korea.}
\email{tkkim@kw.ac.kr}

\author{Hyuck-In Kwon}
\address{Department of Mathematics, Kwangwoon University, Seoul 139-701, Republic of Korea.}
\email{sura@kw.ac.kr}

\author{Jong-Jin Seo}
\address{Department of Mathematics, Pukyong National University, Busan 608-737, Pusan, Republic of Korea.}
\email{seo2011@pknu.ac.kr}

\subjclass{05A10, 05A19.}

\maketitle

\begin{abstract}
In this paper, we consider the degenerate Frobenius-Euler polynomials and investigate some identities of those polynomials.
\end{abstract}

\section{Introduction}

For $u \in{\mathbb{C}}$ with $u \neq 1$, the {\it{Frobenius-Euler polynomials}} are defined by the generating function to be
\begin{equation}\label{1}
\frac{1-u}{e^t-u}e^{xt}=\sum_{n=0} ^{\infty}H_n(x|u)\frac{t^n}{n!},{\text{ (see [1-10])}}.
\end{equation}
When $x=0$, $H_n(x|u)=H_n(0|u)$ are called the {\it{Frobenius-Euler numbers}}.

From \eqref{1}, we have
\begin{equation}\label{2}
H_n(x|u)=\sum_{l=0} ^n \binom{n}{l}H_l(u)x^{n-l},~(n \geq0),{\text{ (see [6-8])}}.
\end{equation}
Note that
\begin{equation*}
\frac{d}{dx}H_n(x|u)=nH_{n-1}(x|u),~(n \in {\mathbb{N}}).
\end{equation*}

In \cite{03}, L. Carlitz define the degenerate Bernoulli polynomials which are given by the generating function to be
\begin{equation}\label{3}
\frac{t}{(1+\lambda t)^{\frac{1}{\lambda}}-1}(1+\lambda t)^{\frac{x}{\lambda}}=\sum_{n=0} ^{\infty}\beta_n(x|\lambda)\frac{t^n}{n!}.
\end{equation}
When $x=0$, $\beta_n(\lambda)=\beta_n(0|\lambda)$ are called the {\it{degenerate Bernoulli numbers}}.
From \eqref{3}, we have
\begin{equation}\label{4}
\lim_{\lambda \rightarrow 0} \frac{t}{(1+\lambda t)^{\frac{1}{\lambda}}-1}(1+\lambda t)^{\frac{x}{\lambda}}=\frac{t}{e^t-1}e^{xt}=\sum_{n=0} ^{\infty}B_n(x)\frac{t^n}{n!},
\end{equation}
where $B_n(x)$ are called {\it{Bernoulli polyomials}}

By \eqref{3} and \eqref{4}, we get
\begin{equation*}
\lim_{\lambda \rightarrow 0}\beta_n(x|\lambda)=B_n(x),~(n \geq 0).
\end{equation*}

The {\it{Stirling number of the first kind}} is defined as
\begin{equation}\label{5}
(x)_n=\sum_{l=0} ^nS_1(n,l)x^l,~(n \geq 0),
\end{equation}
where $(x)_n=x(x-1)\cdots(x-n+1)$, $(x)_0=1$.

By \eqref{3}, we easily get
\begin{equation}\label{6}
\beta_n(x|\lambda)=\sum_{l=0} ^n \binom{n}{l}\beta_{n-l}(x|\lambda)(x|\lambda)_l,~(n \geq 0),
\end{equation}
where $(x|\lambda)_n=x(x-\lambda)\cdots(x-(n-1)\lambda)$.

The {\it{Stirling numbers of the second kind}} is defined by
\begin{equation}\label{7}
x^n=\sum_{l=0} ^n S_2(n,l)(x)_l,~(n\geq 0).
\end{equation}

In this paper, we consider the degenerate Frobenius-Euler polynomials and investigate some properties and identities of those polynomials.

\section{Degenerate Frobenius-Euler polynomials}

For $u \in{\mathbb{c}}$ with $u \neq 1$, we consider degenerate Frobenius-Euler polynomials which are given by the generating function to be
\begin{equation}\label{2-1}
\frac{1-u}{(1+\lambda t)^{\frac{1}{\lambda}}-u}(1+\lambda t)^{\frac{x}{\lambda}}=\sum_{n=0} ^{\infty}h_{n,\lambda}(x|u)\frac{t^n}{n!}.
\end{equation}
When $x=0$, $h_{n,\lambda}(u)=h_{n,\lambda}(0|u)$ are called {\it{degenerate Frobenius-Euler numbers}}.

From \eqref{2-1}, we have
\begin{equation}\label{2-2}
\begin{split}
\frac{1-u}{(1+\lambda t)^{\frac{1}{\lambda}}-u}(1+\lambda t)^{\frac{x}{\lambda}}=&\left(\sum_{l=0} ^{\infty}\frac{h_{l,\lambda}(u)}{l!}t^n\right)\left(\sum_{m=0} ^{\infty}(x|\lambda)_m\frac{t^m}{m!}\right)\\
=&\sum_{n=0} ^{\infty}\left(\sum_{l=0} ^n \binom{n}{l}h_{l,\lambda}(u)(x|\lambda)_{n-l}\right)\frac{t^n}{n!}.
\end{split}
\end{equation}
Thus by \eqref{2-1} and \eqref{2-2}, we get
\begin{equation}\label{2-3}
h_{n,\lambda}(x|u)=\sum_{l=0} ^n\binom{n}{l}h_{l,\lambda}(u)(x|\lambda)_{n-l}.
\end{equation}
From \eqref{2-1}, we can derive the following recurrence relation:
\begin{equation}\label{2-4}
\begin{split}
1-u=&(1+\lambda t)^{\frac{1}{\lambda}}\sum_{n=0} ^{\infty}h_{n,\lambda}(u)\frac{t^n}{n!}-u\sum_{n=0} ^{\infty}h_{n,\lambda}(u)\frac{t^n}{n!}\\
=&\sum_{n=0} ^{\infty}\left(\sum_{l=0} ^n \binom{n}{l} h_{n-l,\lambda}(u)(1|\lambda)_l\right)\frac{t^n}{n!}-u\sum_{n=0} ^{\infty}h_{n,\lambda}(u)\frac{t^n}{n!}\\
=&\sum_{n=0} ^{\infty}\left(h_{n,\lambda}(1|u)-uh_{n,\lambda}(u)\right)\frac{t^n}{n!}.
\end{split}
\end{equation}
By comparing of the coefficients on the both sides of \eqref{2-4}, we get
\begin{equation}\label{2-5}
h_{n,\lambda}(1|u)-uh_{n,\lambda}(u)=
\begin{cases}
1-u, &{\text{if }}n=0\\
0, & {\text{if }} n>0.
\end{cases}
\end{equation}
Note that
\begin{equation}\label{2-6}
\lim_{\lambda\rightarrow 0}\frac{1-u}{(1+\lambda t)^{\frac{1}{\lambda}}-u}(1+\lambda t)^{\frac{x}{\lambda}}=\frac{1-u}{e^t-u}e^{xt}=\sum_{n=0} ^{\infty}H_n(x|u)\frac{t^n}{n!}.
\end{equation}
From \eqref{2-1} and \eqref{2-6}, we have
\begin{equation}\label{2-7}
H_n(x|u)=\lim_{\lambda\rightarrow 0}h_{n,\lambda}(x|u),~(n\geq 0).
\end{equation}
We observe that
\begin{equation}\label{2-8}
\begin{split}
&\frac{1-u}{(1+\lambda t)^{\frac{1}{\lambda}}-u}(1+\lambda t)^{\frac{x+1}{\lambda}}-\frac{u(1-u)}{(1+\lambda t)^{\frac{1}{\lambda}}-u}(1+\lambda t)^{\frac{x}{\lambda}}\\
=&(1-u)(1+\lambda t)^{\frac{x}{\lambda}}=(1-u)\sum_{n=0} ^{\infty}(x|\lambda)_n\frac{t^n}{n!}.
\end{split}
\end{equation}
Therefore, by \eqref{2-1} and \eqref{2-8}, we get
\begin{equation}\label{2-9}
h_{n,\lambda}(x+1|u)-uh_{n,\lambda}(x|u)=(x|\lambda)_n(1-u),~(n\geq0).
\end{equation}
Therefore, we obtain the following theorem.
\begin{theorem}\label{thm1}
For $n \geq 0$, we have
\begin{equation*}
h_{n,\lambda}(x|u)=\sum_{l=0} ^n \binom{n}{l}h_{l,\lambda}(u)(x|\lambda)_{n-l}
\end{equation*}
and
\begin{equation*}
\frac{h_{n,\lambda}(x+1|u)}{1-u}-\frac{u}{1-u}h_{n,\lambda}(x|u)=(x|\lambda)_n.
\end{equation*}
In particular,
\begin{equation*}
h_{n,\lambda}(1|u)-uh_{n,\lambda}(u)=(1-u)\delta_{0,n},
\end{equation*}
where $\delta_{n,k}$ is the Kronecker's symbol.
\end{theorem}
By \eqref{2-1}, we get
\begin{equation}\label{2-10}
\begin{split}
\sum_{n=0} ^{\infty}h_{n,-\lambda}(u)\frac{t^n}{n!}=&\frac{1-u}{(1-\lambda t)^{-\frac{1}{\lambda}}-u}=\frac{1-u}{1-u(1-\lambda t)^{\frac{1}{\lambda}}}(1-\lambda t)^{\frac{1}{\lambda}}\\
=&\frac{1-u^{-1}}{(1-\lambda t)^{\frac{1}{\lambda}}-u^{-1}}(1-\lambda t)^{\frac{1}{\lambda}}=\sum_{n=0} ^{\infty}h_{n,\lambda}(1|u^{-1})(-1)^n\frac{t^n}{n!}.
\end{split}
\end{equation}
In general,
\begin{equation}\label{2-11}
\begin{split}
\sum_{n=0} ^{\infty}h_{n,-\lambda}(-x|u)\frac{t^n}{n!}=&\frac{1-u}{(1-\lambda t)^{-\frac{1}{\lambda}}-u}(1-\lambda t)^{\frac{x}{\lambda}}=\frac{1-u^{-1}}{(1-\lambda t)^{\frac{1}{\lambda}}-u^{-1}}(1-\lambda t)^{\frac{x+1}{\lambda}}\\
=&\sum_{n=0} ^{\infty}h_{n,\lambda}(x+1|u^{-1})(-1)^n\frac{t^n}{n!}.
\end{split}
\end{equation}
Therefore, by \eqref{2-10} and \eqref{2-11}, we obtain the following theorem.
\begin{theorem}\label{thm2}
For $n \geq 0$, we have
\begin{equation*}
(-1)^nh_{n,-\lambda}(-x|u)=h_{n,\lambda}(x+1|u^{-1}).
\end{equation*}
In particular,
\begin{equation*}
(-1)^nh_{n,-\lambda}(u)=(-1)^nh_{n,\lambda}(1|u^{-1}).
\end{equation*}
\end{theorem}

We observe that
\begin{equation}\label{2-12}
\begin{split}
-\frac{1-u}{(1+\lambda t)^{\frac{1}{\lambda}}-u}(1+\lambda t)^{\frac{x}{\lambda}}=&\frac{u^d(1-u^{-1})}{(1+\lambda t)^{\frac{d}{\lambda}}-u^d}\sum_{l=0} ^{d-1}u^{-l}(1+\lambda t)^{\frac{l+x}{\lambda}}\\
=&\frac{u^d(1-u^{-1})}{1-u^d}\sum_{l=0} ^{d-1}\frac{u^{-l}(1-u^d)}{(1+\lambda t)^{\frac{d}{\lambda}}-u^d}(1+\lambda t)^{\frac{l+x}{\lambda}}\\
=&\frac{(1-u^{-1})}{1-u^d}\sum_{l=0} ^{d-1}\sum_{n=0} ^{\infty}u^{d-l}h_{n,\frac{\lambda}{d}}\left(\left.\frac{l+x}{d}\right|u^d\right)d^n\frac{t^n}{n!}\\
=&\sum_{n=0} ^{\infty}\left\{\left(\frac{1-u^{-1}}{1-u^d}\right)d^n\sum_{a=0} ^{d-1}u^{d-a}h_{n,\frac{\lambda}{d}}\left(\left.\frac{a+x}{d}\right|u^d\right)\right\}\frac{t^n}{n!}.
\end{split}
\end{equation}
Therefore, by \eqref{2-1} and \eqref{2-12}, we obtain the following theorem.

\begin{theorem}\label{thm3}
For $n \geq 0$, $d\in{\mathbb{N}}$, we have
\begin{equation*}
h_{n,\lambda}(x|u)=d^n\left(\frac{-1+u^{-1}}{1-u^d}\right)\sum_{a=0} ^{d-1}u^{d-a}h_{n,\frac{\lambda}{d}}\left(\left.\frac{a+x}{d}\right|u^d\right).
\end{equation*}
\end{theorem}
For $r\in{\mathbb{N}}$, let us consider the degenerate Frobenius-Euler polynomials of order $r$ which are given by the generating function to be
\begin{equation}\label{2-13}
\left(\frac{1-u}{(1+\lambda t)^{\frac{1}{\lambda}}-u}\right)^r(1+\lambda t)^{\frac{x}{\lambda}}=\sum_{n=0} ^{\infty}h_{n,\lambda} ^{(r)}(x|u)\frac{t^n}{n!}.
\end{equation}
When $x=0$, $h_{n,\lambda} ^{(r)}(u)=h_{n,\lambda} ^{(r)}(0|u)$ are called {\it{degenerate Frobenius-Euler numbers of order $r$}}.

From \eqref{2-13}, we can easily derive the following equation.
\begin{equation}\label{2-14}
h_{n,\lambda} ^{(r)}(x)=\sum_{l=0} ^n \binom{n}{l}h_{l,\lambda} ^{(r)}(u)(x|\lambda)_{n-l},~(n\geq 0).
\end{equation}
By \eqref{2-13}, we get
\begin{equation}\label{2-15}
\begin{split}
\sum_{n=0} ^{\infty}h_{n,\lambda} ^{(r)}(x+y|\lambda)\frac{t^n}{n!}=&\left(\sum_{l=0} ^{\infty}h_{l,\lambda} ^{(r)}(x)\frac{t^l}{l!}\right)\left(\sum_{m=0} ^{\infty}(x|\lambda)_m\frac{t^m}{m!}\right)\\
=&\sum_{n=0} ^{\infty}\left(\sum_{l=0} ^n h_{l,\lambda} ^{(r)}(x)(x|\lambda)_{n-l}\binom{n}{l}\right)\frac{t^n}{n!}.
\end{split}
\end{equation}
Therefore, by \eqref{2-15}, we obtain the following theorem.
\begin{theorem}\label{thm4}
For $n \geq 0$, we have
\begin{equation*}
h_{n,\lambda} ^{(r)}(x+y)|\lambda)=\sum_{l=0} ^n \binom{n}{l}h_{l,\lambda} ^{(r)}(x)(x|\lambda)_{n-l}.
\end{equation*}
\end{theorem}

From Theorem \ref{thm4}, we note that $h_{n,\lambda} ^{(r)}(x|u)$ is a Sheffer sequence.

By replacing  $t$ by $\frac{1}{\lambda}(e^{\lambda t}-1)$ in \eqref{2-13}, we get
\begin{equation}\label{2-16}
\begin{split}
\left(\frac{1-u}{e^t-u}\right)^{r}e^{xt}=&\sum_{n=0} ^{\infty}h_{n,\lambda} ^{(r)}(x|u)\frac{1}{n!}\left(\frac{1}{\lambda}(e^{\lambda t}-1)\right)^n\\
=&\sum_{n=0} ^{\infty}h_{n,\lambda} ^{(r)}(x|u)\lambda^{-n}\sum_{m=n} ^{\infty}S_2(m,n)\lambda^m\frac{t^m}{m!}\\
=&\sum_{m=0} ^{\infty}\left(\sum_{n=0} ^m h_{n,\lambda} ^{(r)}(x|u)\lambda^{m-n}S_2(m,n)\right)\frac{t^m}{m!}.
\end{split}
\end{equation}

As is well known, the higher-order Frobenius-Euler polynomials are defined by the generating function to be
\begin{equation}\label{2-17}
\left(\frac{1-u}{e^t-u}\right)^re^{xt}=\sum_{n=0} ^{\infty}H_n ^{(r)}(x|u)\frac{t^n}{n!}.
\end{equation}
Thus, by \eqref{2-16} and \eqref{2-17}, we get
\begin{equation}\label{2-18}
H_n ^{(r)}(x|u)=\sum_{n=0} ^m h_{n,\lambda}(x|u)\lambda^{m-n}S_2(m,n).
\end{equation}
Therefore, by \eqref{2-18}, we obtain the following theorem.
\begin{theorem}\label{thm5}
For $m \geq 0$, we have
\begin{equation*}
H_n ^{(r)}(x|u)=\sum_{n=0} ^m h_{n,\lambda}^{(r)}(x|u)\lambda^{m-n}S_2(m,n).
\end{equation*}
\end{theorem}
By replace $t$ by $\log(1+\lambda t)^{\frac{1}{\lambda}}$ in \eqref{2-17}, we have
\begin{equation}\label{2-19}
\begin{split}
\left(\frac{1-u}{(1+\lambda t)^{\frac{1}{\lambda}}-u}\right)^r(1+\lambda t)^{\frac{x}{\lambda}}=&\sum_{n=0} ^{\infty}H_n ^{(r)}(x|u)\frac{1}{n!}\frac{1}{\lambda^n}\left(\log(1+\lambda t)\right)^n\\
=&\sum_{n=0} ^{\infty}H_n ^{(r)}(x|u)\lambda^{-n}\sum_{m=n} ^{\infty}S_1(m,n)\lambda^m\frac{t^n}{n!}\\
=&\sum_{m=0} ^{\infty}\left(\sum_{n=0} ^m H_n ^{(r)}(x|u)\lambda^{m-n}S_1(m,n)\right)\frac{t^m}{m!}.
\end{split}
\end{equation}
Therefore, by \eqref{2-13} and \eqref{2-18}, we obtain the following theorem.
\begin{theorem}\label{thm6}
For $m \geq 0$, we have
\begin{equation*}
h_{n,\lambda} ^{(r)}(x|u)=\sum_{n=0} ^m H_n ^{(r)}(x|u)\lambda^{m-n}S_1(m,n).
\end{equation*}
\end{theorem}
We observe that
\begin{equation}\label{19}
\begin{split}
&\left(\frac{1-u}{(1+\lambda t)^{\frac{1}{\lambda}}-u}\right)^r(1+\lambda t)^{\frac{x+1}{\lambda}}-u\left(\frac{1-u}{(1+\lambda t)^{\frac{1}{\lambda}}-u}\right)^r(1+\lambda t)^{\frac{x}{\lambda}}\\
=&\left(\frac{1-u}{(1+\lambda t)^{\frac{1}{\lambda}}-u}\right)^{r-1}(1+\lambda t)^{\frac{x}{\lambda}}(1-u)=\sum_{n=0} ^{\infty}h_{n,\lambda} ^{(r-1)}(x|u)(1-u)\frac{t^n}{n!}.
\end{split}
\end{equation}
Therefore, by \eqref{2-13}, we obtain the following theorem.
\begin{theorem}\label{thm7}
For $n \geq 0$, we have
\begin{equation*}
\frac{1}{1-u}\left\{h_{n,\lambda} ^{(r)}(x+1|u)-uh_{n,\lambda} ^{(r)}(x|u)\right\}=h_{n,\lambda} ^{(r-1)}(x|u).
\end{equation*}
\end{theorem}

\noindent{\bf{{\scshape Remark. }}} Let $u=-1$. Then, by \eqref{2-1}, we get
\begin{equation}\label{20}
\frac{1}{t}\frac{2t}{(1+\lambda t)^{\frac{1}{\lambda}}+1}(1+\lambda t)^{\frac{x}{\lambda}}=\sum_{n=0} ^{\infty}h_{n,\lambda}(x|-1)\frac{t^n}{n!}.
\end{equation}
Thus, by \eqref{20}, we get
\begin{equation*}
\frac{2t}{(1+\lambda t)^{\frac{1}{\lambda}}+1}(1+\lambda t)^{\frac{x}{\lambda}}=t\sum_{n=0} ^{\infty}h_{n,\lambda}(x|-1)\frac{t^n}{n!}.
\end{equation*}
Now, we define the degenerate Genocchi polynomials which are given by the generating function to be
\begin{equation}\label{21}
\frac{2t}{(1+\lambda t)^{\frac{1}{\lambda}}+1}(1+\lambda t)^{\frac{x}{\lambda}}=\sum_{n=0} ^{\infty}g_{n,\lambda}(x)\frac{t^n}{n!}.
\end{equation}
From \eqref{20} and \eqref{21}, we have
\begin{equation*}
g_{0,\lambda}(x)=0{\text{ and }}h_{n,\lambda}(x|-1)=\frac{1}{n+1}g_{n+1,\lambda}(x),~(n\geq 0).
\end{equation*}

\end{document}